\documentclass[11pt]{amsart}
\usepackage{amssymb,amscd}

\newcommand{\cA}{{\mathcal A}}

\newcommand{\cE}{{\mathcal E}}

\newcommand{\cV}{{\mathcal V}}
\newcommand{\cW}{{\mathcal W}}
\newcommand{\cZ}{{\mathcal Z}}
\newcommand{\fB}{{\mathfrak B}}
\newcommand{\fT}{{\mathfrak T}}
\newcommand{\fU}{{\mathfrak U}}
\newcommand{\fb}{{\mathfrak b}}
\newcommand{\fg}{{\mathfrak g}}

\newcommand{\ft}{{\mathfrak t}}
\newcommand{\fu}{{\mathfrak u}}
\newcommand{\mA}{{\mathbb A}}
\newcommand{\mC}{{\mathbb C}}

\newcommand{\mP}{{\mathbb P}}

\newcommand{\p}{{\Phi}}
\newcommand{\w}{{\Omega}}

\newcommand{\A}{{\alpha}}
\newcommand{\B}{{\beta}}

\newcommand{\z}{{\zeta}}

\newcommand{\Ad}{\operatorname{Ad}}
\newcommand{\End}{\operatorname{End}}
\newcommand{\Lie}{\operatorname{Lie}}
\newcommand{\Tr}{\operatorname{Tr}}

\newtheorem{theorem}{Theorem}
\newtheorem{corollary}{Corollary}
\newtheorem{proposition}{Proposition}
\newtheorem{lemma}{Lemma}

\newtheorem{Remark}{Remark}
\def\thebricarbibliography#1{\par\bigskip

\begin{center}
{\normalsize \bf References}
\end{center}
\par
\noindent\list
 {[\arabic{enumi}]}{\settowidth\labelwidth{[#1]}\leftmargin\labelwidth
 \advance\leftmargin\labelsep
 \usecounter{enumi}}

 \sloppy\clubpenalty4000\widowpenalty4000
 \sfcode`\.=1000\relax}

\title{The equivariant cohomology ring of regular varieties}
\author{Michel~Brion}
\address{Universit\'e de Grenoble I\\
D\'epartement de Math\'ematiques\\
Institut Fourier, UMR 5582 du CNRS\\
38402 Saint--Martin d'H\`eres Cedex, France}
\email{Michel.Brion@ujf-grenoble.fr}
\author{James~B.~Carrell}
\thanks{{The second author was partially supported by the Natural 
Sciences and Engineering Research Council of Canada}}
\address{Department of Mathematics\\
The University of British Columbia\\
Vancouver, B.C. \\
Canada, V6T 1Z2}
\email{carrell@math.ubc.ca}

\begin{document}

\maketitle

\begin{abstract} 
Let $\fB$ denote the upper triangular subgroup of $SL_2(\mC)$, 
$\fT$ its diagonal torus and $\fU$ its unipotent radical. A complex
projective variety $Y$ endowed with an algebraic action of $\fB$ 
such that the fixed point set $Y^{\fU}$ is a single point, is called
regular. Associated to any regular $\fB$--variety $Y$, there is a
remarkable affine curve $\cZ_Y$ with a $\fT$--action which was 
studied in \cite{CRELLE}. In this note, we show that the coordinate
ring $\mC[\cZ_Y]$ is isomorphic with the equivariant cohomology ring
$H_{\fT}^*(Y)$ with complex coefficients, when $Y$ is smooth or, more
generally, is a $\fB$--stable subvariety of a regular smooth
$\fB$--variety $X$ such that the restriction map from $H^*(X)$ to
$H^*(Y)$ is surjective. This isomorphism is obtained as a refinement
of the localization theorem in equivariant cohomology; it applies e.g.
to Schubert varieties in flag varieties, and to the Peterson variety
studied in \cite{kost}. Another application of our isomorphism is a
natural algebraic formula for the equivariant push forward.
\end{abstract}

\section{Preliminaries} 
Let $\fB$ be the group of upper triangular $2\times 2$ complex matrices
of determinant $1$. Let $\fT$ (resp.~$\fU$) be the subgroup of $\fB$
consisting of diagonal (resp.~unipotent) matrices. We have
isomorphisms $\lambda:\mC^*\rightarrow \fT$ and  
$\varphi:\mC\rightarrow \fU$, where
$$
\lambda(t)=\begin{pmatrix}t &0\\0 & t^{-1} \end{pmatrix}
$$ 
and 
$$
\varphi(u)=\left(\begin{matrix}
1&u\cr 0&1\cr
\end{matrix}\right),
$$
together satisfying the relation  
\begin{eqnarray}\label{GMGA}
\lambda(t)\varphi(u)\lambda(t^{-1})=\varphi(t^2u).
\end{eqnarray}
Consider the  generators 
$$\cV=\dot{\varphi}(0)=\begin{pmatrix}
0&1\\ 0&0\end{pmatrix} \quad \text{and} 
\quad \cW=\dot{\lambda}(1)=\begin{pmatrix}1 & 0\\0&-1\end{pmatrix}$$ 
of the Lie algebras $\Lie(\fU)$ and $\Lie(\fT)$ respectively.
Then $[\cW,\cV]=2\cV$, and
\begin{eqnarray}\label{AD} \Ad(\varphi(u))\cW=\cW-2u\cV
\end{eqnarray}
for all $u\in\mC$.

In this note, $X$ will denote a smooth complex projective algebraic
variety endowed with an algebraic action of $\fB$ such that the fixed
point scheme $X^{\fU}$ consists of one point $o$. Both the
$\fB$--variety $X$ and the action will be called {\em regular}. Then
$o\in X^{\fT}$, since $X^{\fU}$ is $\fT$--stable. Moreover, $X^{\fT}$
is finite by Lemma 1 of \cite{CRELLE}. Thus we may write
\begin{eqnarray}\label{TFPS}
X^{\fT}=\{\z_1=o,\z_2,\ldots,\z_r\}.
\end{eqnarray}
Clearly $r=\chi(X)$, the Euler characteristic of $X$. 

Let $H^*_{\fT}(X)$ denote the $\fT$--equivariant cohomology ring of
$X$ with complex coefficients. To define it, let $\cE$ be a
contractible space with a free action of $\fT$ and let 
$X_{\fT}=(X\times \cE)/\fT$ 
(quotient by the diagonal $\fT$--action). Then 
$$
H_{\fT}^*(X)=H^*(X_{\fT}).
$$
It is well known that the equivariant cohomology ring $H^*_{\fT}(pt)$ of
a point is the polynomial ring $\mC[z]$, where $z$ denotes the linear
form on the Lie algebra of $\fT$ such that $z(\cW)=1$. The degree of $z$
is $2$. Thus $H_{\fT}^*(X)$ is a graded algebra over the polynomial
ring $\mC[z]=H^*_{\fT}(pt)$ (via the constant map $X\to pt$). 

In our situation, the restriction map in cohomology
$$
i_{\fT}^*:H^*_{\fT}(X)\to H^*_{\fT}(X^{\fT}),
$$
induced by the inclusion $i:X^{\fT}\hookrightarrow X$,
is injective \cite{GKM}. By (\ref{TFPS}), 
$$
H^*_{\fT}(X^{\fT})= \bigoplus_{j=1}^r H^*_{\fT}(\z_j)\cong
\bigoplus_{j=1}^r \mC[z],
$$
so each $\alpha \in H_{\fT}^*(X)$ defines an $r$--tuple of polynomials
$(\alpha_{\z_1},\ldots,\alpha_{\z_r})$. That is,
\begin{eqnarray}\label{REST}
i_{\fT}^*(\alpha)=(\alpha_{\z_1},\ldots,\alpha_{\z_r}).
\end{eqnarray}
We will define a refined version of this restriction map in \S 4 below.

\section{The $\fB$--stable curves} 

Throughout this note, a curve in $X$ will be a purely one--dimensional
closed subset of $X$. One which is stable under a subgroup $G$ of $\fB$
is called a $G$--curve. The $\fB$--curves in $X$ play a crucial role, 
so we will next establish a few of their basic properties.

\begin{proposition}\label{prop1} 
If $X$ is a regular $\fB$--variety, then every irreducible
$\fB$--curve $C$ in $X$ has the form $C=\overline{\fB \cdot \z_j}$ for
some index $j\ge 2$. Moreover, every $\fB$--curve contains $o$. In
particular, there are only finitely many irreducible $\fB$--curves in
$X$, and they all meet at $o$. 
\end{proposition}

\begin{proof} It is clear that if $j\ge 2$, then 
$C=\overline{\fB \cdot \z_j}$ is a $\fB$--curve in $X$ containing $\z_j$,
which, by the Borel Fixed  Point Theorem, also contains $o$, since $o$
is the only $\fB$--fixed point. Conversely, every $\fB$--curve $C$ in
$X$ contains $o$, and at least one other $\fT$--fixed point. 
Indeed, since $o$ has an affine open $\fT$--stable neighbourhood
$X_o$ in $X$, the complement $C - X_o$  is nonempty and $\fT$--stable.
It follows immediately that $C=\overline{\fB \cdot x}$ for some
$x\in X^{\fT} - o$.
\end{proof}

Consider the action of $\fB$ on the projective line ${\mP}^1$ given by
$$
\left(\begin{matrix}t&u\cr 0&t^{-1}\cr\end{matrix}\right)\cdot z
= \frac{z}{t(t+uz)}
$$
(the inverse of the standard action). It has $(\mP^1)^{\fT}=\{0,\infty\}$
and $(\mP^1)^{\fU}=\{0\}$ and is therefore regular. Note that if
$u\in\mC^*$, 
$$
\left(\begin{matrix}t&u\cr 0&t^{-1}\cr\end{matrix}\right)\cdot \infty
= \frac{1}{tu},
$$
so $\varphi(u)\cdot \infty=u^{-1}$.

The diagonal action of $\fB$ on $X\times \mP^1$ is also regular. 
By Proposition 1, the irreducible $\fB$--curves in $X\times \mP^1$ are of
the form $\overline{\fB\cdot (x,\infty)}$ or $\overline{\fB\cdot (x,0)}$,
where $x\in X^{\fT}$. Only the first type will play a role here. Thus put
$Z_j=\overline{\fB \cdot (\z_j,\infty)}$, and let $\pi_j:Z_j\to\mP^1$ be
the second projection.  Clearly every $\pi_j$ is bijective, hence $Z_j
\cong \mP^1$. In addition,
$$
Z_j=\{(\varphi(u)\cdot \z_j,u^{-1})\mid u\in
\mC^*\}\cup \{(\z_j,\infty)\}\cup \{(o,0)\},
$$ 
so $Z_i\cap Z_j=\{(o,0)\}$ as long as $i\ne j$. 
Moreover, restricting $\pi_j$ gives an isomorphism
$$
p_j:Z_j - (\z_j,\infty)\to\mA^1.
$$ 
Finally, put 
$$
Z=\bigcup _{1\le j\le r}Z_j.
$$
Thus $Z$ is the union of all irreducible $\fB$--stable curves 
in $X\times\mP^1$ that are mapped onto $\mP^1$ by the second
projection $\pi:X\times\mP^1\rightarrow \mP^1$.

\bigskip
\section{The fundamental scheme $\cZ$} 
Let $\cA$ denote the vector field on $X\times \mA^1$ defined by
\begin{eqnarray}\label{DEFA} \cA _{(x,v)}=2\cV_x-v\cW_x. 
\end{eqnarray}
Obviously $\cA$ is tangent to the fibres of the projection to $\mA^1$.
By (\ref{AD}),
\begin{eqnarray}\label{AID} (\Ad(\varphi(u))\cW)_x=-u\cA_{(x,u^{-1})}.
\end{eqnarray}
The contraction operator $i(\cA)$ defines a sheaf of ideals
$i(\cA)(\Omega_{X\times \mA^1}^1)$
of the structure sheaf of $X\times \mA^1$. Let
$\cZ$ denote the associated closed subscheme of $X\times\mA^1$. In
other words, $\cZ$ is the zero scheme of $\cA$. 

\begin{Remark} The vector field $\cA$ here is a variant of the vector
field studied in \cite{CRELLE}. Both have the same zero scheme.
\end{Remark}

The properties of $\cZ$ figure prominently in this note. First put 
\begin{eqnarray}\label{BMC}
X_o=\{x\in X \mid \lim_{t\to \infty} \lambda(t)\cdot x = o\}.
\end{eqnarray}
Clearly, $X_o$ is $\fT$--stable; and it follows easily from (\ref{GMGA})
that $X_o$ is open in $X$ (see Proposition 1 of \cite{CRELLE} for
details). Hence by the Bialynicki--Birula decomposition theorem
\cite{BB}, $X_o$ is $\fT$--equivariantly isomorphic to the tangent space
$T_oX$, where $\fT$ acts by its canonical representation at
a fixed point. The weights of the associated action of $\lambda$ on
$T_oX$ are all negative. So we may choose coordinates $x_1,\dots,x_n$
on $X_o\cong T_o X$ that are eigenvectors of $\fT$; the weight $a_i$ of
$x_i$ is a positive integer (it turns out to be even, 
see \cite{PNAS}). This identifies the positively graded ring
$\mC[X_o]$ with $\mC[x_1,\ldots,x_n]$, where $\deg x_i = a_i$. Now
$X\times\mP^1$ contains $X_o\times\mA^1$ as a $\fT$--stable affine open
subset, with coordinate ring $\mC[x_1,\ldots,x_n,v]$, where $v$ has
degree $2$.

\begin{proposition}\label{zeroscheme} 
The scheme $\cZ$ is reduced, and contained in $X_o\times \mA^1$ as a
$\fT$--curve. Its ideal in $\mC[X_o\times\mA^1]=\mC[x_1,\ldots,x_n,v]$
is generated by 
\begin{eqnarray}\label{IZ} 
v\cW(x_1)-2\cV(x_1),\dots,v\cW(x_n)-2\cV(x_n).
\end{eqnarray} 
These form a homogeneous regular sequence in
$\mC[x_1,\ldots,x_n,v]$, and the degree of each $v\cW(x_i)-2\cV(x_i)$
equals $a_i+2$. The irreducible components of $\cZ$ are the $\cZ_j$,
$1\le j \le r$, where 
\begin{eqnarray}\label{IRREDCOMP}
\cZ_j= Z_j - \{(\z_j,\infty)\}. 
\end{eqnarray}
Each $\cZ_j$ is mapped isomorphically to $\mA^1$ by the second
projection $p$. In particular, $p$ is finite and flat of degree
$r$, and $\cZ$ has $r$ irreducible components.
Any two such components meet only at $(o,0)$. 
\end{proposition}

\begin{proof}
This follows from the results in \S 3 of \cite{CRELLE}. We provide
direct arguments for the reader's convenience. 

Since $[\cW,\cV]=2\cV$ and $\cW(v)=2v$ (on $\mA^1$), it follows that
$\cA$ commutes with the vector field induced by the diagonal
$\fT$--action on $X\times\mA^1$. Thus, $\cZ$ is $\fT$--stable. 
 
Next we claim that (\ref{IRREDCOMP}) holds set--theoretically. Let 
$(x,v)\in X\times\mA^1$ with $v\ne 0$. Then $(x,v)\in\cZ$ if and only
if $\cW_x-2v^{-1}\cV_x=0$, that is, $x$ is a zero of 
${\rm Ad}(\varphi(v^{-1})\cW)$; equivalently, 
$\varphi(-v^{-1})\cdot x\in X^{\fT}$. 
On the other hand, $(x,0)\in\cZ$ if and only if $\cV_x=0$, that is,
$x=o$. Thus $\cZ= Z \cap(X\times\mA^1)$ (as sets). Further,
$Z\cap(X\times\mA^1)$ is an open affine $\fT$--stable neighborhood of
$(o,0)$ in $Z$, and hence equals $Z\cap(X_o\times\mA^1)$. This implies
our claim.

It follows that $\cZ\subseteq X_o\times\mA^1$ (as schemes), so that the
ideal of $\cZ$ is generated by
$v\cW(x_1)-2\cV(x_1),\dots,v\cW(x_n)-2\cV(x_n)$. 
These polynomials are homogeneous of degrees $a_1+2,\dots,a_n+2$;
together with $v$, they have only the origin as their common zero
(since $o$ is the unique zero of $\cV$). Hence
$$
v\cW(x_1)-2\cV(x_1),\dots,v\cW(x_n)-2\cV(x_n),v
$$ 
form a regular sequence in $\mC[x_1,\dots,x_n,v]$, and $v$ is a
non--zero divisor in $\mC[\cZ]$. As a consequence, the $\mC[v]$--module
$\mC[\cZ]$ is finitely generated and free. In other words, $p:\cZ\to
\mA^1$ is finite and flat.

Fix $v_0\ne 0$ and consider the scheme--theoretic intersection
$\cZ\cap(X\times v_0)$. This identifies with the zero scheme of 
$2\cV_x - v_0\cW_x$ in $X$, that is, to the zero scheme of 
${\rm Ad}(\varphi(v_0^{-1})\cW)$. The latter consists of $r=\chi(X)$
distinct points, so that $\cZ\cap(X\times v_0)$ is reduced. Since
$\mC[\cZ]$ is a free module over $\mC[v]$, it follows easily that 
$\cZ \cap(v\ne 0)$ is reduced. But the open subset $\cZ \cap(v\ne 0)$
is dense in $\cZ$ by (\ref{IRREDCOMP}), and $\cZ$ is a complete
intersection in $\mA^{n+1}$. Thus, $\cZ$ is reduced; this completes
the proof.
\end{proof}

\bigskip
\section{The refined restriction}  
We now define our refined restriction on equivariant cohomology. 
Let $\alpha\in H^*_{\fT}(X)$. Recall from (\ref{REST}) that 
$i^*_{\fT}(\alpha)=(\alpha_{\z_1},\dots,\alpha_{\z_r})$,
where each $\alpha_{\z_j}\in \mC[z]$. We regard each $\alpha_{\z_j}$
as a polynomial function on $\cZ_j$ (isomorphic to $\mA^1$ via $p$),
and hence on $\cZ_j - (o,0)$. Since $\cZ - (o,0)$ is the disjoint
union of the $\cZ_j - (o,0)$, this yields an algebra homomorphism 
$$
\rho:H^*_{\fT}(X)\to\mC[\cZ - (o,0)]
$$
such that $\rho(\alpha)(x,v)=\alpha_{\z_j}(v)$ whenever 
$(x,v)\in\cZ_j - (o,0)$. In particular, $\rho(z)(x,v)=v$, so that
$\rho(z)=v$. And since $i^*_{\fT}$ preserves the grading, the
same holds for $\rho$. 

Note that the value $\alpha_{\z_j}(0)$ at the origin is independent of
the index $j$. (For $X_{\fT}=(X\times \cE)/\fT$ is connected since
both $X$ and $\cE$ are, so $H^0(X_{\fT})$ is by definition the set of
constant functions on $X_{\fT}$. Consequently, the component of
$\alpha$ in degree 0 gives the same value at each fixed point.) 
Now let $\mC_0[\cZ]$ denote the subalgebra of $\mC[\cZ - (0,o)]$
consisting of all elements that extend continuously to $\cZ$
in the classical topology. Then
$$
\rho\big(H_{\fT}^*(X)\big)\subseteq \mC_0[\cZ].
$$
We will show below that $\rho\big(H_{\fT}^*(X)\big)$ is in fact the
coordinate ring $\mC[\cZ]$. To do this, we compute the image
under $\rho$ of an equivariant Chern class.

\section{Equivariant Chern theory} 
Recall that if $Y$ is an algebraic variety with an action of an
algebraic group $G$, a vector bundle $E$ on $Y$ is said to be
$G$--linearized if there is an action of $G$ on $E$ lifting that on
$Y$, such that each $g\in G$ defines a linear map from $E_y$ to
$E_{g\cdot y}$, for any $y\in Y$. In particular, if $y\in Y^G$, then
we have a representation of $G$ in $E_y$, and hence a representation
of $\Lie(G)$. Thus, each $\xi\in \Lie(G)$ acts on $E_y$, by
$\xi_y\in\End(E_y)$.

Also, recall that the $k$th equivariant Chern class 
$c_k^G(E)\in H_G^{2k}(Y)$ is defined to be the $k$th Chern class of
the vector bundle  
$$
E_G =(E\times \cE)/G \to X_G = (X\times \cE)/G,
$$
where $\cE$ is a contractible space with a free action of $G$. 
For $y\in Y^G$, the restriction $c_k^G(E)_y$ lies in
$H^*_G(pt)$. The latter identifies with a subring of the coordinate
ring of $\Lie(G)$, and we have
$$
c_k^G(E)_y(\xi)=\Tr_{\wedge^k E_y}(\xi_y)
$$
for any $\xi\in\Lie(G)$.

Returning to our previous situation, let $E$ be a $\fB$--linearized
vector bundle on $X$. Then each $(x,v)\in \cZ$ is a zero of 
$v\cW-2\cV\in\Lie(\fB)$. This yields an element 
$(v\cW -2\cV)_x \in \End(E_x)$ 

\begin{lemma}\label{chern}
Let $E$ be a $\fB$--linearized vector bundle on $X$ and let $k$ be a
non--negative integer. Then we have for any $(x,v)\in\cZ$:
\begin{eqnarray}\label{RHO} 
\rho(c_k^T(E))(x,v)= \Tr_{\wedge^k E_x}((v\cW -2\cV)_x).
\end{eqnarray}
As a consequence, $\rho(c_k^T(E))\in \mC[\cZ]$.
\end{lemma}

\begin{proof}
It suffices to check (\ref{RHO}) for $v\neq 0$. Let $j$ be the index 
such that $x=\varphi(v^{-1})\cdot \z_j$. Letting 
$\cW_{\z_j}\in \End(E_{\z_j})$ denote the lift of $\cW$  at
$\z_j$,  we have 
\begin{eqnarray}\label{RHO'}
\rho(c_k^T(E))(x,v)=c_k^T(E)_{\z_j}(v)=
v^k\Tr_{\wedge^k E_{\z_j}}(\cW_{\z_j}).
\end{eqnarray}
Since $E$ is $\fB$--linearized, this equals
$$\displaylines{
v^k \Tr_{\wedge^k E_x}({\rm Ad}(\varphi(v^{-1}))\cW)_x)=
v^k \Tr_{\wedge^k E_x}((\cW-2v^{-1}\cV)_x)
\hfill\cr\hfill
=\Tr_{\wedge^k E_x}((v\cW-2\cV)_x),
\cr}$$
which proves (\ref{RHO}).
\end{proof}

We now obtain some of the main properties of $\rho$.

\begin{proposition}\label{comp}
The image of the morphism $\rho:H^*_{\fT}(X)\to\mC[\cZ - (o,0)]$  
is contained in $\mC[\cZ]$.
\end{proposition}

\begin{proof} Since $X$ is smooth and projective, we have
an exact sequence 
\begin{eqnarray}\label{EXSEQ} 
0\rightarrow zH^*_{\fT}(X)\rightarrow H_{\fT}^*(X)\rightarrow 
H^*(X)\rightarrow 0.
\end{eqnarray}
By Proposition 3 of \cite{CRELLE}, the algebra $H^*(X)$ is
generated by Chern classes of $\fB$--linearized vector bundles on $X$,
so by Nakayama's Lemma, $H^*_{\fT}(X)$ is generated (as an algebra) by
their equivariant Chern classes. It follows from this and Lemma 1 that
$\rho\big(H_{\fT}^*(X)\big)\subseteq \mC[\cZ]$.
\end{proof}

\bigskip
\section{ The first main result} 
We now prove

\begin{theorem}\label{main}
For a smooth projective regular variety $X$, the homomorphism  
$$
\rho:H^*_{\fT}(X)\rightarrow \mC[\cZ]
$$ 
is an isomorphism of graded algebras.
\end{theorem}

\begin{proof} 
To see $\rho$ is injective, suppose $\rho(\alpha)=0$.
Then each $\alpha_{\z_j}$ is $0$, so $i_{\fT}^*(\alpha)=0$, where
$i_{\fT}^*$ is the restriction (cf. (\ref{REST})).
But, as noted previously,  $i_{\fT}^*$ is injective, so $\alpha=0$.

To complete the proof, it suffices to check that the Poincar\'e series
of $H^*_{\fT}(X)$ and $\mC[\cZ]$ coincide (since both algebras are
positively graded, and $\rho$ preserves the grading). Denote the
former by $F_X(t)$ and the latter by $F_{\cZ}(t)$. By (\ref{EXSEQ}),
we have an isomorphism 
\begin{eqnarray}\label{ESC} H^*_{\fT}(X)/zH^*_{\fT}(X)\cong H^*(X).
\end{eqnarray} 
Since $H^*_{\fT}(X)$ is a free $\mC[z]$--module and $z$ has degree $2$,
(\ref{ESC}) implies 
\begin{eqnarray}\label{EQFOR}
F_X(t)=\frac{P_X(t)}{1-t^2}
\end{eqnarray} 
where $P_X(t)$ is the Poincar\'e polynomial of $H^*(X)$.
On the other hand, since $p:\cZ\to\mA^1$ is finite, flat and
$\fT$--equivariant (where $\fT$ acts linearly on $\mA^1$ with weight $2$),
we have
$$
F_{\cZ}(t)=\frac{P_{\cZ}(t)}{1-t^2}
$$
where $P_{\cZ}(t)$ is the Poincar\'e polynomial of the
finite--dimensional graded algebra $\mC[\cZ]/(v)$.
We now use the fact that the cohomology 
ring of $X$ is isomorphic as a graded algebra to the coordinate
ring of the zero scheme of $\cV$ with the principal grading
\cite{CRELLE}. So
\begin{eqnarray}\label{AVF}
H^*(X)\cong \mC[x_1,\dots,x_n]/(\cV(x_1),\dots,\cV(x_n))
\cong \mC[\cZ]/(v),
\end{eqnarray}
where the second isomorphism follows from Proposition 2. Thus, 
$P_X(t)=P_{\cZ}(t)$, whence $F_X(t)=F_{\cZ}(t)$.
\end{proof}

For a simple example, let $X=\mP^n$. Let $e:\mC^{n+1} \to \mC^{n+1}$
denote the nilpotent linear transformation defined by 
$$
v_n\to v_{n-1}\to \dots \to v_1\to v_0 \to 0,
$$
where $(v_0,v_1,\ldots, v_n)$ is the standard basis of $\mC^{n+1}$.
Also, let $h=\text{diag}(n+1,n-1, \dots ,-n+1, -n-1)$. Then
$[h,e]=2e$, so we obtain a $\fB$--action on $\mP^n$. 
This action is regular with unique fixed point $[v_0]=[1,0,\dots,0]$;
its neighborhood $X_o$ is the standard affine chart centered at $[v_0]$.
Let $(x_1, \dots ,x_n)$ be the usual affine coordinates at
$[v_0]$. That is, $x_j=z_j/z_0$, where $[z_0,z_1, \dots,z_n]$
are the homogeneous coordinates on $\mP^n$. Then each $x_j$ is
homogeneous of degree $2j$, and a straightforward computation gives 
\begin{equation*} 
e(x_1)=x_2 -x_1^2, \, e(x_2)=x_3 -x_1 x_2,  
\dots , \, e(x_{n-1})=x_n-x_1 x_{n-1}, 
\end{equation*}
$$ e(x_n)= -x_1x_n.$$
Thus, the ideal of $\cZ$ in $\mC[x_1,\ldots,x_n,v]$ is generated by 
$-x_2+x_1(x_1+v)$, $-x_3+x_2(x_1+2v)$, $\ldots$,
$-x_n+x_{n-1}(x_1+(n-1)v)$, $x_n(x_1+nv)$. By Theorem 1, it follows that
$$
H_{\fT}^*(\mP ^n) \cong \mC[x_1,v]/\big(\prod _{m=0}^n (x_1 + mv)\big).
$$
And by equivariant Chern theory, $x_1=-c_1^{\fT}(L)$, 
where $L$ is the tautological line bundle on $\mP^n$.
This presentation of $H_{\fT}^*(\mP ^n)$ is well known, and can be
derived directly.

\bigskip
\section{The general case}
In this section, we will make some comments on regular
$\fB$--varieties, dropping the smoothness assumption. Let $Y$ be a
complex projective variety endowed with a $\fB$--action such that
$X^{\fU}$ is a unique point $o$. If $Y$ is singular, it is necessarily
singular at $o$, and so the $\fT$--stable neighbourhood $Y_o$ of $o$
defined in (\ref{BMC}) is singular. Hence the results in \S 3 on the
structure of $\cZ$ (relative to $Y$) do not necessarily obtain. To be
able to conclude something, let us assume that $Y$ is equivariantly
embedded into a smooth projective regular $\fB$--variety $X$. Thus the
curve $\cZ$ (relative to $X$) is well defined and enjoys all the
properties derived above. We can therefore define
\begin{eqnarray}\label{ZSUBY}
\cZ _Y=\cZ\cap (Y_o\times \mA^1)
\end{eqnarray}
taking the intersection to be reduced. In other words, $\cZ_Y$ is 
the union of the components of $\cZ$ lying in $Y_o\times \mA^1$. Now the
construction of \S 4 yields a graded homomorphism
$$
\rho_Y:H_{\fT}^*(Y)\rightarrow \mC_0[\cZ_Y].
$$

If we make the additional assumption that $H^*(Y)$ is generated by
Chern classes of $\fB$--linearized vector bundles, then the odd
cohomology of $Y$ is trivial, so that $H^*_{\fT}(Y)$ is a free module
over $\mC[z]$, and the restriction
$i^*_{\fT}:H_{\fT}^*(Y)\to H_{\fT}^*(Y^{\fT})$ is
injective \cite{GKM}. Hence, $\rho_Y$ is injective. And by exactly the
same argument, $\rho_Y(H_{\fT}^*(Y)) \subseteq \mC[\cZ_Y]$.
The obstruction to $\rho_Y$ being an isomorphism is therefore
the equality of Poincar\'e series of these graded algebras. 

Now there is a natural situation where this assumption
is satisfied, hence one can generalize Theorem \ref{main}.
Assume that the inclusion $j:Y\to X$ induces a surjection
$j^*:H^*(X)\to H^*(Y)$. Then $H^*(Y)$ is generated by Chern classes of
$\fB$--linearized vector bundles, and from the Leray spectral
sequence, one sees that $j^*_{\fT}:H_{\fT}^*(X)\to H_{\fT}^*(Y)$ is
surjective. The surjectivity of $j^*$ holds, for example, in the case
of Schubert varieties in flag varieties. Also, see \S 8.

Proceeding as above, we obtain an extension of Theorem 1.

\begin{theorem}\label{gen-main} 

Suppose $X$ is a smooth projective variety with a regular $\fB$--action
and $Y$ is a closed $\fB$--stable subvariety for which the restriction
map $H^*(X)\rightarrow H^*(Y)$ is surjective. 
Then the map $\rho_Y:H^*_{\fT}(Y)\rightarrow \mC_0[\cZ_Y]$
yields a graded algebra isomorphism
$$
\rho_Y:H^*_{\fT}(Y)\rightarrow \mC[\cZ_Y]
$$
fitting into a commutative diagram
\begin{eqnarray}\label{CD1}
\begin{CD} H_{\fT}^*(X) @>{\rho}>> \mC[\cZ] \\
                  @VVV           @VVV  \\
              H_{\fT}^*(Y)    @>{\rho_Y}>> \mC[\cZ_Y]
\end{CD}
\end{eqnarray}
where the vertical maps are the natural restrictions. Both vertical
maps are surjections.

\end{theorem}

\begin{proof}
We know $\rho$ is an isomorphism from Theorem 1, and $\rho_Y$ is
injective by the above remarks. Since the restriction map 
$\mC[\cZ] \rightarrow \mC[\cZ_Y]$ is clearly surjective, it therefore
follows  that $\rho_Y$ is also. 
\end{proof}

As a corollary, we obtain one of the main results of \cite{CRELLE}.

\begin{corollary} 
With the notation and assumptions of Theorem 2, there
exists a commutative diagram
\begin{eqnarray}\label{CD2}
\begin{CD} \mC[\cZ]/(v) @>{\psi}>> H^*(X) \\
                  @VVV           @VVV  \\
              \mC[\cZ_Y]/(v)    @>{\psi_Y}>> H^*(Y)
\end{CD}
\end{eqnarray}
where $\psi$ and $\psi_Y$ are graded algebra isomorphisms
and the vertical maps are the natural restrictions. 
\end{corollary}

Note that $\mC[\cZ_Y]/(v)$ is the coordinate ring of the schematic
intersection of $\cZ_Y$ and $X\times 0$ in $X\times\mA^1$.


\bigskip

\section{Equivariant cohomology of the Peterson variety}
Let $G$ be a complex semi--simple linear algebraic group. Fix a pair of
opposite Borel subgroups $B$ and $B^-$ and let $T=B\cap B^-$, a
maximal torus of $G$. Denote the corresponding Lie algebras by $\fg$,
$\fb$, $\fb ^-$ and $\ft$, and let $\p^+$ and $\p^-$ be the
roots of the pair $(G,T)$ which arise from $\fb$ and $\fb^-$ respectively. 

Let $M$ be a $B$-submodule of $\fg$ containing $\fb$. Then 
$$
M=\fb\oplus\bigoplus_{\alpha\in \Omega(M)} \fg_{\alpha}, \text{ where }
\w(M)=\{ \A\in \p^- \mid \fg _{\A}\subset M\}.
$$
So $\w(M)$ is the set of weights of the quotient $M/\fb$.
The $B$--module $M\subseteq \fg$ yields a homogeneous vector bundle
$G\times^B M$ over $G/B$, together with a morphism 
$G\times^B M\rightarrow \fg$ induced by $(g,m)\mapsto gm$. 
The fiber of this morphism at an arbitrary
$x\in\fg$ identifies with 
$$
Y_M(x)=\{gB\in G/B \mid g^{-1}x\in M\}.
$$
This is a closed subvariety of $G/B$, stable under the action of 
$G_{\mC x}$ (the isotropy group of the line $\mC x$).

If $x$ is regular and semi--simple, then $Y_M(x)$ is known as a 
{\em Hessenberg variety}. It was shown in \cite{proc} that $Y_M(x)$ is 
nonsingular, and its Poincar\'e polynomial was also determined.

On the other hand, if $x$ is a regular nilpotent, the situation is
quite different. To describe it, recall that any regular nilpotent is
conjugate to a principal nilpotent, defined as follows. Let
$\A_1,\dots,\A_\ell \in \p^+$ denote the simple roots, and fix a
non--zero $e_{\A_j}$ in each $T$--weight space
$\fg_{\A_j}\subset\fb$. Then a principal nilpotent is an  element of
the form $e=\sum_{j=1}^\ell e_{\A_j}$. By choosing an element
$h\in\ft$ for which $\A_j(h)=2$ for each index $j$, we obtain a pair
$(e,h) \in \fb \times \ft$ 
which determines a regular action of $\fB$ on $G/B$ having unique
$\fB$--fixed point $o=B$. Since $\fB$ fixes the line $\mC e$, it
stabilizes the subvariety $Y_M(e)$ for any $M$. 
Note that the semi--simple element $h$ is also regular.

Let $\fT\subset \fB$ denote the maximal torus in $\fB$ with Lie
algebra $\mC h$. Then $Y_M(e)^{\fT} =(G/B)^T \cap Y_M(e)$. Let 
$W=N_G(T)/T$ be the Weyl group of $(G,T)$, and choose a representative
$n_w\in N_G(T)$ for each $w\in W$. Then 
$(G/B)^T =\{n_w B~\vert~w\in W\}$, and $n_wB\in Y_M(e)^{\fT}$ if and
only if $n_w^{-1} e\in M$. This is equivalent to the condition
$w^{-1}(\A_j)\in \w (M)\cup \p^+$ for each $j$. 

The cohomology algebras of the varieties $Y_M(e)$ have been determined
by Dale Peterson. Formulated in our terms, his result goes as follows.

\begin{theorem} 
For any $B$--submodule $M$ of $\fg$ containing $\fb$, the restriction map 
$$
H^*(G/B) \to H^*(Y_M(e))
$$ 
is surjective. Hence $($by Theorem 2$)$ 
$$
H_{\fT}^*(Y_M(e))\cong \mC[\cZ_{Y_M(e)}].
$$
Moreover, $\cZ_{Y_M(e)}$ is a complete intersection. Thus,
the cohomology ring $H^*(Y_M(e))$ satisfies Poincar\'e duality. 
Its Poincar\'e polynomial is given by the product formula
\begin{equation}\label{PF} P(Y_M(e),t)= \prod _{-\A\in \w(M)} 
\frac{1-t^{\text{ht}(\A) +1}}{1-t^{\text{ht}(\A)}},
\end{equation}
where $ht(\A)$ denotes the sum of the coefficients of $\A\in\Phi^+$
over the simple roots. 
\end{theorem}

In the case where $M=\fg$, (\ref{PF}) is a well known product formula 
for the Poincar\'e polynomial of the flag variety.

The ideals $I(\cZ_{Y_M (e)})$ admit explicit expressions
related to the geometry of the nilpotent variety of $\fg$.
Let $U^-\subset G$ denote the unipotent radical of $B^-$. Since the
natural map $\mu:U^-\to G/B$, $\mu(u)=uB$, is $T$--equivariant with
respect to the conjugation action of $T$ on $U^-$, we may
equivariantly identify $U^-$ and the open cell $(G/B)_o$. 
Let $\fu^-$ denote the Lie algebra of $U^-$,  and let 
$\Pi_*:\fg \to \fu^-$ denote the projection. 
To get this explicit picture, we first note that 
$L_{u^{-1}*}\cV_u=\Pi_*\big(u^{-1}e\big)$, where $L_u$ denotes left
translation by $u\in U^-$, and $L_{u*}$ is its differential
(cf. \cite[\S 2.1]{JDG}). Consequently, 
\begin{equation}\label{}L_{u^{-1}*}\cA_{(u,v)}
=2\Pi_*\big(u^{-1}e\big) -v L_{u^{-1}*}\cW_u.
\end{equation}
Defining  $F_\alpha (u)$ as the component 
of $L_{u^{-1}*}\cA_{(u,v)}$ in $\fg_\alpha$, it follows that
$$
I(\cZ)=(F_\alpha ~\vert~ \alpha \in \Phi^-).
$$
More precisely, we can write
$$
u^{-1}e= e +k(u) +\sum _{\alpha\in\Phi^-} v_\alpha (u) e_\A,
$$
where $k\in \ft\otimes\mC[U^-]$ and all $v_\alpha \in \mC[U^-]$. 
Then $F_\A=2v_\A -vw_\A$, where we have put 
$w_\A (u)=(L_{u^{-1}*}\cW_u)_\A.$ Hence,
$$
H_\fT ^*(G/B,\mC)\cong 
\mC[U^- \times \mA^1]/(2v_\A -vw_\A ~\vert~ \alpha \in \Phi^-).
$$

Since the functions $v_\B$ ($\B\in \Phi^- -\w(M)$)
cut out $Y_M(e)_o\times \mA^1$, it follows
$I(\cZ_{Y_M (e)})$ is the ideal of the variety 
cut out by the $2v_\A -vw_\A$, where $ \A\in \Phi^-$, and the
$v_\B$, where $\B\in \Phi^- -\w(M).$

The case where 
\begin{equation}\label{PM} M=
\fb\oplus\bigoplus _{\A\in\Phi^+,~\A(h)>2} \fg _{-\A}.
\end{equation}
is of particular interest. Then, as shown by Kostant in
\cite[cf. (9)]{kost}, the Giventhal--Kim and Peterson formulas for the
flag variety quantum cohomology  
\cite{GK,P} may be interpreted as asserting an isomorphism of graded
rings 
$$
\mC[Y_M(e)_o]\cong QH^*(G^\vee/B^\vee),
$$
where $G^\vee$ and $B^\vee$ denote the Langlands duals of $G$ and $B$
respectively and $QH^*(G^\vee/B^\vee)$ is the complex quantum
cohomology ring of $G^\vee/B^\vee$. This led Kostant to call $Y_M(e)$
the {\em Peterson variety}. Notice that by Theorem 3,
$P(Y_M(e),t)=(1+t)^{n-1}$.

\bigskip
\section{An alternative proof of Theorem 1}
We now give another proof of Theorem \ref{main} independent of
(\ref{AVF}). Hence we will also obtain another proof of
(\ref{AVF}). Denote by $[X^{\fT}]_{\fT}\in H^*_{\fT}(X)$ the 
equivariant cohomology class of the $\fT$--fixed point set. We compute
its image under $\rho:H^*_{\fT}(X)\to\mC[\cZ]$.

\begin{lemma}\label{class}
$\rho([X^{\fT}]_{\fT})$ is the restriction to $\cZ$ of the Jacobian
determinant of the polynomial functions
$v\cW(x_1)-2\cV(x_1),\dots,v\cW(x_n)-2\cV(x_n)$ in the variables
$x_1,\dots,x_n$.
\end{lemma}

\begin{proof}
Let $T_X$ be the tangent bundle to $X$, then $\cW$ yields a
$\fT$--invariant global section of $T_X$, with zero scheme $X^{\fT}$. 
Thus, we have $c_n^{\fT}(T_X)=[X^{\fT}]_{\fT}$ in $H^*_{\fT}(X)$. 
On the other hand, $T_X$ carries a $\fB$--linearization,
and we have by Lemma \ref{chern}:
$$
\rho(c_n^{\fT}(T_X))(x,v)=
\Tr_{\wedge^n T_x X}(v\cW_x-2\cV_x).
$$
This is the Jacobian determinant of
$v\cW(x_1)-2\cV(x_1),\dots,v\cW(x_n)-2\cV(x_n)$.
\end{proof}

We will also need the following easy result of commutative algebra; we
will provide a proof, for lack of a reference.

\begin{lemma}\label{jac}
Let $P_1,\dots,P_n$ be polynomial functions in $x_1,\dots,x_n$, that
are weighted homogeneous for a positive grading defined by 
$\deg x_j = a_j$. If the origin is the unique common zero to
$P_1,\ldots,P_n$, then the Jacobian determinant $J(P_1,\ldots,P_n)$ is
not in the ideal $(P_1,\ldots,P_n)$.
\end{lemma}

\begin{proof}
We begin with the special case where $P_1,\ldots,P_n$ are
homogeneous. We argue by induction on $n$, the result being
evident for $n=1$. Let $d_1,\ldots,d_n$ be the degrees of
$P_1,\ldots,P_n$. We have the Euler identities
$$\displaylines{
x_1\frac{\partial P_1}{\partial x_1}+\cdots +
x_n\frac{\partial P_1}{\partial x_n}= d_1P_1,
\cr\cdots\cr
x_1\frac{\partial P_n}{\partial x_1}+\cdots +
x_n\frac{\partial P_n}{\partial x_n}=d_nP_n
\cr}$$
that we view as a system of linear equalities in $x_1,\ldots,x_n$. Its
determinant is the Jacobian $J(P_1,\ldots,P_n)=J$. For 
$1\leq i,j\leq n$, let $J_{i,j}$ be its maximal minor associated with
the $i$th line and the $j$th column. Then we have
\begin{eqnarray}\label{Jac}
x_i J=\sum_{j=1}^n (-1)^{j-1}d_j P_j J_{ij}.
\end{eqnarray}
Assume that $J\in (P_1,\ldots,P_n)$ and write 
$$
J=f_1P_1+\cdots+f_nP_n
$$ 
with $f_1,\ldots,f_n\in\mC[x_1,\ldots,x_n]$. Using \ref{Jac} for
$i=1$, it follows that $(x_1f_1-d_1J_{11})P_1$ is in the ideal
$(P_2,\ldots,P_n)$. But $P_1,\ldots,P_n$ form a regular sequence in
$\mC[x_1,\ldots,x_n]$, as they are homogeneous and the origin is their
unique common zero. Therefore, 
$x_1f_1-d_1J_{11}\in (P_2,\ldots,P_n)$. In other words, 
$J_{11}\in (x_1,P_2,\ldots,P_n)$.

After a linear change of coordinates, we can assume that
$(x_1,P_2,\ldots,P_n)$ is a regular sequence. For $2\leq i\leq n$, let
$$
Q_i(x_2,\ldots,x_n)=P_i(0,x_2,\ldots,x_n).
$$ 
Then we have in $\mC[x_2,\ldots,x_n]$:
$$
J_{11}(0,x_2,\cdots,x_n)\in (Q_2,\ldots,Q_n).
$$
Now $J_{11}(0,x_2,\ldots,x_n)=J(Q_2,\ldots,Q_n)$, and $Q_2,\ldots,Q_n$
are homogeneous polynomial functions of $x_2,\ldots,x_n$ having the
origin as their unique common zero. But this contradicts the inductive 
assumption. This completes the proof in the homogeneous case.

Consider now the case where $P_1,\ldots,P_n$ are quasi--homogeneous for
the weights $a_1,\ldots,a_n$. Let $y_1,\ldots,y_n$ be indeterminates;
then the functions 
$$
(y_1,\ldots,y_n)\mapsto P_i(y_1^{a_1},\ldots,y_n^{a_n})
$$
$( 1\leq i\leq n)$ are homogeneous polynomials with the origin as
their unique commmon zero. Their Jacobian determinant is
$$
(\prod_{i=1}^n a_i y^{a_i-1})\; 
J(P_1,\ldots,P_n)(y_1^{a_1},\ldots,y_n^{a_n}).
$$
By the first step of the proof, this function of $(y_1,\ldots,y_n)$ is
not in the ideal generated by the
$P_i(y_1^{a_1},\ldots,y_n^{a_n})$. Thus,
$J(P_1,\ldots,P_n)(x_1,\ldots,x_n)$ cannot be in $(P_1,\ldots,P_n)$.
\end{proof}

\noindent
{\it The proof of Theorem 1.}
The functions $\cV(x_1),\ldots, \cV(x_n)$ satisfy the assumption of Lemma
\ref{jac}, since they are quasi--homogeneous for the grading defined by
the action of $\fT$, and $o$ is the unique zero of $\cV$. Thus,
the Jacobian determinant of these functions is not in the ideal that
they generate. By Lemma \ref{class}, it follows that
$\rho([X^{\fT}]_{\fT})$ is not divisible by $v$ in $\mC[\cZ]$, for
$\mC[\cZ]/(v)=\mC[x_1,\ldots,x_n]/(\cV(x_1),\ldots,\cV(x_n))$.

The $\mC[z]$--linear map $\rho:H^*_{\fT}(X)\to\mC[\cZ]$ defines a map 
$$
\overline{\rho}:H^*(X)\cong H^*_{\fT}(X)/(z)\to\mC[\cZ]/(v),
$$
a graded ring homomorphism. It suffices to prove that
$\overline{\rho}$ is an isomorphism. Note that the spaces $H^*(X)$ and
$\mC[\cZ]/(v)$ have dimension $r$, so that it suffices to check
injectivity of $\overline{\rho}$.

The image in $H^*(X)$ of $[X^{\fT}]_{\fT}$ is $r[pt]$, where $[pt]$ denotes
the cohomology class of a point. Thus, $\overline{\rho}([pt])\neq 0$. 
Let $(\alpha_1,\ldots,\alpha_r)$ be a basis of $H^*(X)$ consisting of
homogeneous elements and let $(\beta_1,\ldots,\beta_r)$ be the dual
basis for the intersection pairing 
$(\alpha,\beta)\mapsto\int_X (\alpha\cup\beta)\cap[X]$.
Then the homogeneous component of degree $2n$ in each product
$\alpha_i\cup\beta_j$ equals $[pt]$ if $i=j$, and is zero
otherwise. Assume that 
$\overline{\rho}(t_1\alpha_1+\cdots+t_r\alpha_r)=0$ for some complex
numbers $t_1,\ldots,t_r$. Multiplying by $\overline{\rho}(\beta_j)$
and taking the homogeneous component of degree $2n$, we obtain
$t_j\overline{\rho}([pt])=0$ whence $t_j=0$. Thus, $\overline{\rho}$
is injective, and the proof is complete.
\medskip

\bigskip
\section{The equivariant push forward}
Next we describe the equivariant push--forward map
$$
\int_X:H^*_{\fT}(X)\to H^*_{\fT}(pt)=\mC[z],
$$
associated to the map $X\to pt$. 
Note that $\int_X$ is $\mC[z]$--linear and homogeneous of degree $-2n$.
Denote by $J$ the restriction to $\cZ$ of the Jacobian determinant of
the polynomial functions 
$$
v\cW(x_1)-2\cV(x_1),\ldots,v\cW(x_n)-2\cV(x_n)
$$ 
in the variables $x_1,\ldots,x_n$. Then $J$ is homogeneous of degree
$2n$, as follows from Lemma 2.

\begin{theorem}\label{residues}
For any $f\in\mC[\cZ]$, the function
$$
v\mapsto \sum_{(x,v)\in\cZ} \frac{f(x,v)}{J(x,v)}
$$
is polynomial. Further, for any $\alpha\in H^*_{\fT}(X)$:
\begin{eqnarray}\label{EPF}
(\int_X \alpha)(v)=\sum_{(x,v)\in\cZ}\frac{\rho(\alpha)(x,v)}{J(x,v)}.
\end{eqnarray}
\end{theorem}

\begin{proof}
Since $\mC[\cZ]$ is graded free $\mC[v]$--module of rank
$r$, we may choose a homogeneous basis $f_1,\ldots,f_r$ with $f_1=1$. 
As $J\notin(v)$, we may also assume $f_r=J$. If $f\in \mC[\cZ]$, let 
$\varphi(f)$ be the $r$th coordinate of $f$ in this basis.
Then the map
$$
\mC[\cZ]\times\mC[\cZ]\to \mC[v],~(f,g)\mapsto \varphi(fg)
$$
is a non--degenerate bilinear form, since it reduces modulo $(v)$ to
the duality pairing on 
$\mC[\cZ]/(v)=\mC[x_1,\ldots,x_n]/(\cV(x_1),\ldots,\cV(x_n))$. 
If $(g_1,\ldots,g_r)$ is the dual basis with respect to
$(f_1,\ldots,f_r)$ for this bilinear form, then $g_1,\ldots,g_r$ are
homogeneous and satisfy 
${\rm deg}(f_i)+{\rm deg}(g_i)=2n$ for all $i$. As a consequence,
the kernel of $\varphi$ is generated by $f_1,\ldots,f_{r-1}$ and also
by $g_2,\ldots,g_r$ as a $\mC[z]$--module. Since $J=f_r$, we have
$\varphi(Jg_1)=\cdots=\varphi(Jg_{r-1})=0$, so that
$\varphi(Jf_2)=\cdots=\varphi(Jf_r)=0$ whereas
$\varphi(Jf_1)=\varphi(J)=1$.

Let 
$$
\Tr:\mC[\cZ]\to\mC[v]
$$ 
be the trace map for the (finite, flat) morphism $p:\cZ\to\mA^1$. Then
$$
\Tr(f)(v)=\sum_{(x,v)\in\cZ} f(x,v)
$$ 
for all $v\in\mA^1$; in particular, $\Tr(1)=r$. Since 
$\Tr$ is homogeneous of degree $0$ and $\mC[z]$--linear, its
kernel is a graded complement of the $\mC[v]$--module
$\mC[v]=\mC[v]f_1$ in $\mC[\cZ]$. It follows that this kernel is
generated by $f_2,\ldots,f_r$. Thus, we have
$$
\Tr(g)=r\varphi(Jg)
$$
for all $g\in\mC[\cZ]$. This equality holds then for all
$g\in\mC[\cZ][v^{-1}]$, and hence for all rational functions on $\cZ$. 
Note that $J$ restricts to a non--zero function on any component of
$\cZ$, so that $\frac{1}{J}$ is a rational function on $\cZ$. 
We thus have
$$
\Tr(\frac{f}{J})= r\varphi(f)
$$
for any $f\in\mC[\cZ]$; this implies the first assertion. The second
assertion follows from the localization theorem in equivariant
cohomology.
\end{proof}

\section{Some concluding remarks}
If $\alpha\in H^*_{\fT}(X)$ is a product of equivariant Chern classes,
then (\ref{EPF}) is the equivariant Bott residue formula for regular
$\fB$--varieties (cf. \cite{BOTT}).

One can also, by similar methods extend (\ref{EPF}) to obtain a
formula for the equivariant Gysin homomorphism associated to an
equivariant morphism of regular smooth projective $B$--varieties
(cf. \cite{IJM} for the case of flag varieties).

More generally, consider a smooth projective variety $X$ with an
action of an arbitrary torus $T$ such that $X^T$ is finite. Then the
precise version of the localization theorem given in \cite{GKM} yields
a reduced affine scheme whose coordinate ring is the equivariant 
cohomology ring of $X$ (see \cite{BRI} for details and applications). 
The scheme $\cZ$ of this note gives a more explicit picture, but the
requirement of a regular $\fB$--action is harder to satisfy. It would be
nice to be able to relax the regularity assumption to allow $X^{\fU}$
to have positive dimension.

\begin{thebricarbibliography}{}

\bibitem{IJM} E.\ Akyildiz  and J. B. \ Carrell: 
{\it An algebraic formula for the Gysin homomorphism from $G/B$ to $G/P$}, 
Ill. J. Math. {\bf 31} (1987), 312--320.

\bibitem{PNAS}  E.\ Akyildiz  and J. B. \ Carrell: 
{\it  A generalization of the Kostant--Macdonald identity}, 
Proc. Nat. Acad. Sci. U.S.A. {\bf 86} (1989), 3934--3937.

\bibitem{BB}  A.\ Bialynicki--Birula:  
{\em Some theorems on actions of algebraic groups}, 
Ann. of Math. {\bf 98} (1973), 480--497.

\bibitem{BOTT} R.\ Bott: 
{\em Vector fields and characteristic numbers}, 
Mich. Math. J. {\bf 14} (1967), 231--244.

\bibitem{BRI} M.\ Brion: 
{\em Poincar\'e duality and equivariant (co)homology.
Dedicated to William Fulton on the occasion of his 60th birthday}, 
Michigan Math. J. {\bf 48} (2000), 77--92.

\bibitem{JDG} J. B.\ Carrell: 
{\em Bruhat cells in the nilpotent variety and the intersection rings
of Schubert varieties}, J. Differential Geometry {\bf 37} (1993),
651--668. 

\bibitem{CRELLE} J. B.\ Carrell: 
{\em Deformation of the nilpotent zero scheme and the intersection 
rings of invariant subvarieties}, 
J. Reine Angew. Math. {\bf 460} (1995), 37--54.

\bibitem{proc} F. \ DeMari, C. \ Procesi and M. Shayman: 
{\em Hessenberg varieties},
Trans. Amer. Math. Soc. {\bf 332} (1992), 529--534.

\bibitem{GK} A. \ Giventhal and B. \ Kim: 
{\em Quantum cohomology of flag manifolds and Toda lattices}, Comm.
Math. Phys. {\bf 168} (1995), 609--641.

\bibitem{GKM} R. M.\ Goresky, R. \ Kottwitz and R. \ MacPherson: 
{\em Equivariant cohomology, Koszul duality and the localization theorem}, 
Invent. Math. {\bf 131} (1998), 25--84.

\bibitem{kost} B. \ Kostant: 
{\em Flag manifold quantum cohomology, the Toda lattice and the
representation with highest weight}, Selecta Mathematica, New Series 
{\bf 2} (1996), 43--91.

\bibitem{P} D. \ Peterson: 
{\em Quantum cohomology of $G/P$}, in preparation.

\end{thebricarbibliography}{}

\bigskip
\noindent
Acknowledgement: The second author thanks Shrawan Kumar and Dale Peterson.
\end{document}